\documentclass[final]{amsart}%{{{1
\usepackage{amsfonts,amssymb,enumerate,xspace}

\usepackage[ps2pdf]{hyperref}
\usepackage[notcite,notref]{showkeys}

\newcommand{\snseqns}{\eqref{e:X-def}--\eqref{e:u-def}\xspace}

\newcommand{\del}{\partial}
\newcommand{\lap}{\triangle}
\newcommand{\inv}{^{-1}}
\newcommand{\transpose}{^\text{t}}
\newcommand{\leqs}{\leqslant}
\newcommand{\geqs}{\geqslant}
\newcommand{\grad}{\nabla}
\newcommand{\gradt}{\grad \transpose}
\newcommand{\divergence}{\grad \cdot}

\newcommand{\W}{\boldsymbol{\mathrm{W}}}
\newcommand{\Wns}{\mathcal{W}}
\newcommand{\E}{\boldsymbol{\mathrm{E}}}
\newcommand{\lhp}{\boldsymbol{\mathrm{P}}}
\newcommand{\D}[2][u]{\del_t #2 + \left(#1 \cdot \grad\right) #2 - \nu \lap #2}

\newcommand{\I}{\mathcal{I}}
\newcommand{\Ifn}{I}
\newcommand{\Imatrix}{\mathbb{I}}
\newcommand{\R}{\mathbb{R}}
\newcommand{\holderspace}[2]{\ensuremath{C^{\ifx0#1{#2}\else{#1,#2}\fi}}}
\newcommand{\cspace}[1]{\ensuremath{C^{#1}}}
\newcommand{\lballT}[3]{\ensuremath{\mathcal{L}^{#1,#2}_{#3}}}
\newcommand{\lball}{\ensuremath{\mathcal{L}}}
\newcommand{\uball}{\ensuremath{\mathcal{U}}}

\newcommand{\ubound}{U}

\newif\iftextstyle
\textstyletrue
\everydisplay\expandafter{\the\everydisplay\textstylefalse}
\newcommand{\norm}[1]{\iftextstyle\|#1\|\else\left\|#1\right\|\fi}
\newcommand{\hnorm}[3]{\iftextstyle\|#1\|\else\left\|#1\right\|\fi_{\ifx0#2{#3}\else{#2,#3}\fi}}
\newcommand{\hsnorm}[2]{\iftextstyle|#1|\else\left|#1\right|\fi_{#2}}
\newcommand{\cnorm}[2]{\iftextstyle\|#1\|\else\left\|#1\right\|\fi_{C^{#2}}}

\numberwithin{equation}{section}
\allowdisplaybreaks

\newtheorem{theorem}{Theorem}[section]
\newtheorem{lemma}[theorem]{Lemma}
\newtheorem{prop}[theorem]{Proposition}
\newtheorem{cor}[theorem]{Corollary}

\newtheorem*{theorem*}{Theorem}
\newtheorem*{prop*}{Proposition}

\theoremstyle{definition}
\newtheorem{defn}[theorem]{Definition}

\theoremstyle{remark}
\newtheorem{rem}[theorem]{Remark}
\newtheorem*{rem*}{Remark}
\begin{document}
%{{{1 Title + etc
\title{A stochastic perturbation of inviscid flows}
\author{Gautam Iyer}
\address{%
Department of Mathematics\\
The University of Chicago\\
Chicago, Illinois 60637}
\email{gautam@math.uchicago.edu}
\begin{abstract}
We prove existence and regularity of the stochastic flows used in the stochastic Lagrangian formulation of the incompressible Navier-Stokes equations (with periodic boundary conditions), and consequently obtain a $\holderspace{k}{\alpha}$ local existence result for the Navier-Stokes equations. Our estimates are independent of viscosity, allowing us to consider the inviscid limit. We show that as $\nu \to 0$, solutions of the stochastic Lagrangian formulation (with periodic boundary conditions) converge to solutions of the Euler equations at the rate of $O(\sqrt{\nu t})$.
%
%We consider a stochastic flow with drift $u$ and diffusion coefficient $\sqrt{2 \nu}$. We demand that the drift be recovered from the flow map using the Weber formula, as in the Eulerian-Lagrangian formulation of the Euler equations. This has been shoIn the absence of diffusion, this will yield the Euler equations. We first prove the existence of such stochastic flows, and that the expected value of this process approximates the Navier-Stokes equations (with viscosity $\nu$) to order $O(t^{3/2})$. As a result of our estimates we also obtain a local existence and uniqueness results for the Navier-Stokes equations.
\end{abstract}
\maketitle

\section{Introduction}%{{{1
Consider an incompressible inviscid fluid with velocity field $u$ in the absence of external forcing. The evolution of the velocity field is governed by the Euler \cite{cm} equations
\begin{gather}
\label{e:euler-evol} \del_t u + (u \cdot \grad) u + \grad p = 0\\
\label{e:incompressible} \divergence u = 0.
\end{gather}
Viscosity introduces a diffusive term in the Euler equations and equation \eqref{e:euler-evol} becomes
\begin{equation}
\label{e:ns-evol} \D{u} + \grad p = 0.
\end{equation}

The Kolmogorov backward and Feynman-Kac formulae \cite{ks} show that any linear, diffusive, second order PDE can be obtained by averaging out a stochastic perturbation of an ODE. The theory for non-linear PDE's is not as well developed. We are interested in interpreting the Navier-Stokes equations as the average of a suitable stochastic perturbation of the Euler equations.

Many interesting non-linear PDE's have been interpreted as averaging of stochastic processes, the simplest example being the Kolmogorov reaction diffusion equation \cite{legall}. In two dimensions the same is possible for the Navier-Stokes equations as the vorticity satisfies a standard Fokker-Plank equation. This combined with the Biot-Savart law led to the random vortex method \cite{mb} and has been used and studied extensively. In three dimensions the problem is a little harder as the vorticity equation is no longer of Fokker-Plank type, and the non-linearity causes trouble.

In \cite{sznit} Le Jan and Sznitman used a backward in time branching process in Fourier space to express the Navier-Stokes equations as the expected value of a stochastic process. This approach led to a new existence theorem, and was later \cite{ossiander} generalized and physical space analogues were developed.

An approach more along the lines of this paper was developed by Busnello, Flandoli and Romito \cite{bf} who considered `noisy' flow paths, and used Girsanov transformations to recover the velocity field. They obtained the $3$-dimensional Navier-Stokes equations in this form, and generalized their method to work for a general class of second order parabolic equations. A different technique was used by Gomes \cite{gomes} to express the diffusive Lagrangian \cite{elns} as the expected value minimizer of a suitable functional. Finally we mention similar systems have been considered by Jourdain et al in \cite{jourdain}.\smallskip

Our approach%
\footnote{In the original version of this paper, our intention was to propose this as a physically meaningful \textit{model} for the Navier-Stokes equations. We presented a proof that the solution of the system considered here differs from the solution of the Navier-Stokes equations by $O(t^{3/2})$. Six months after submission of the original version of this paper, the author and Peter Constantin~\cite{detsns} discovered that the equations considered here are exactly equivalent to the Navier-Stokes equations.}
is to introduce a Brownian drift into the active vector formulation \cite{ele} of the Euler equations. Peter Constantin and the author showed \cite{detsns} that this provides a physically meaningful, explicit stochastic representation of the Navier-Stokes equations. While long time dynamics of the system we consider are presently unknown, we hope that techniques used here will lead to control of the growth of certain quantities with non-zero probability. For example, we would like to find an exponential bound for $\grad X$ which holds with non-zero probability. Finding an almost sure bound of this form will lead to global existence.

% We prove a local existence theorem for the process we consider, and show that this process approximates the Navier-Stokes equations to order $O(t^{3/2})$. The proof relies on the Eulerian-Lagrangian formulation of the Euler \cite{ele} and Navier-Stokes \cite{elns} equations developed by P. Constantin. At present we are unable to find explicit deterministic equations for the averaged process, and are not even sure if such equations exist.

In this paper, we consider the flow given by the stochastic differential equation
\begin{equation}
\label{e:X-def} dX = u \,dt + \sqrt{2 \nu} \,dB
\end{equation}
with initial data
\begin{equation}
\label{e:X-idata} X(a,0) = a.
\end{equation}
Here $\nu >0$ represents the viscosity, and $B$ represents a $3$-dimensional Wiener process (we use the letter $B$ to avoid confusion with the Weber operator). We recover the velocity field from $X$ by
\begin{gather}
\label{e:A-def} A = X\inv\\
\label{e:u-def} u = \E \lhp\left[ (\gradt A) \, (u_0 \circ A) \right]
\end{gather}
where $\E$ denotes the expected value with respect to the Wiener measure, $\lhp$ denotes the Leray-Hodge projection \cite{cm} on divergence free vector fields, and $u_0$ is the deterministic initial data. We clarify that by $X\inv$ in equation \eqref{e:A-def} we mean the spatial inverse of $X$. We impose periodic boundary conditions, though all theorems proved here will also work if we work with the domain $\R^3$ and impose a decay at infinity condition instead.

The motivation for considering the above system arises from the fact that in the absence of viscosity, the system \snseqns reduces to
\begin{gather}
\label{e:ele-adef} \del_t A + (u \cdot \grad) A = 0\\
A(x,0) = x\\
\label{e:ele-webber} u = \lhp\left[ (\gradt A) (u_0 \circ A) \right].
\end{gather}
Peter Constantin proved \cite{ele} that $u$ is a solution of the (deterministic) system \eqref{e:ele-adef}--\eqref{e:ele-webber} if and only if $u$ is a solution of the incompressible Euler equations \eqref{e:euler-evol}--\eqref{e:incompressible} with initial data $u_0$. Thus the system \snseqns can be thought of as superimposing the Wiener process on the flow map, intuitively representing Brownian motion of fluid particles. Physically, the Brownian particle interaction is regarded as the source of viscosity, and the equivalence of \snseqns and the Navier-Stokes equations proved in \cite{detsns} confirms this. We remark that equation \eqref{e:u-def} provides an explicit formula for $u$ in terms of the map $X$.

In this paper, we provide a self contained proof of a $\holderspace{k}{\alpha}$ local existence theorem for the stochastic system \snseqns. The proof in~\cite{detsns} showing equivalence between Navier-Stokes and \snseqns relies crucially on spatial regularity of solutions as stated in Theorem~\ref{t:snsexist}. We remark that the stochastic representation of Busnello, Flandoli and Romito does not admit a self contained existence proof as we have here.

The estimates, and existence time can be chosen independent of the viscosity, thus enabling us to consider the vanishing viscosity limit. We show that as $\nu \to 0$, the solution of \snseqns converges to the solution of the Euler equations at the rate of $O(\sqrt{\nu t})$. We remark that the limit $\nu \to 0$ is not well understood in bounded domains using classical methods. We hope that this stochastic formulation (when extended to bounded domains) will give us a better handle on computing this limit.\smallskip

In the next section, we establish our notational convention, and describe precisely the results we prove in this paper. In section~\ref{s:wbounds} we prove bounds on the Weber operator, which are essential to all proofs presented in this paper. In section~\ref{s:snsexist} we prove local existence for \snseqns and the vanishing viscosity limit. Finally, in section~\ref{s:nsexist}, we digress and present an alternate proof of local existence for the Navier-Stokes equations using the diffusive Lagrangian formulation~\cite{elns}.
%
% Finally we show that the velocity $u$ above approximates the velocity of a viscous fluid (with viscosity $\nu$) to order $O(t^{3/2})$. As a consequence of our estimates we also obtain a \holderspace{k}{\alpha} local existence and uniqueness result for the Navier-Stokes equations. We are at present unaware if the $O(t^{3/2})$ estimate is sharp, or if there is a closed deterministic system of equations the velocity $u$ above satisfies.
%
\section{Notational convention and Description of results}\label{s:dresults}%{{{1

In this section we describe the main results we prove. We begin by establishing our notational convention. We let $\I$ denote the cube $[0,L]^3$ with side of length $L$. We define the H\"older norms and semi-norms on $\I$ by
\begin{gather*}
\hsnorm{u}{\alpha} = \sup_{x,y \in \I} L^\alpha \frac{|u(x) - u(y)|}{|x-y|^\alpha}\\
\cnorm{u}{k} = \sum_{|m|\leqs k} L^{|m|} \sup_\I |D^m u|\\
\hnorm{u}{k}{\alpha} = \cnorm{u}{k} + \sum_{|m| = k} L^k \hsnorm{D^m u}{\alpha}
\end{gather*}
where $D^m$ denotes the derivative with respect to the multi index $m$. We let \cspace{k} denote the space of all $k$-times continuously differentiable spatially periodic functions on $\I$, and \holderspace{k}{\alpha} denote the space of all spatially periodic $k+\alpha$ H\"older continuous functions. The spaces \cspace{k} and \holderspace{k}{\alpha} are endowed with the norms $\cnorm{\cdot}{k}$ and $\hnorm{\cdot}{k}{\alpha}$ respectively.

We use $\Ifn$ to denote the identity function on $\R^3$ or $\I$ (depending on the context), and use $\Imatrix$ to denote the identity matrix. The first theorem we prove addresses local (in time) existence for the system \snseqns:

\begin{theorem}\label{t:snsexist}%{{{2 stochastic local existence
Let $k \geqs 1$ and $u_0 \in \holderspace{k+1}{\alpha}$ be divergence free. There exists a time $T = T(k, \alpha, L, \hnorm{u_0}{k+1}{\alpha})$, but independent of viscosity, and a pair of functions $\lambda, u \in C([0,T], \holderspace{k+1}{\alpha})$ such that $u$ and $X = \Ifn + \lambda$ satisfy the system \snseqns. Further $\exists \ubound = \ubound(k, \alpha, L, \hnorm{u_0}{k+1}{\alpha})$ such that $t \in [0,T] \implies \hnorm{u(t)}{k+1}{\alpha} \leqs \ubound$.
\end{theorem}%}}}

We prove this theorem in section \ref{s:snsexist}. Our proof will also give a local existence result for the Euler equations, or any stochastic perturbation similar to the one considered here. We remark that the estimates required for this theorem along with Constantin's diffusive Lagrangian formulation \cite{elns} also gives us local existence for the Navier-Stokes equations. In section \ref{s:nsexist}, we digress and present this proof.

We remark that Theorem \ref{t:snsexist} is still true when $k=0$. The only modification we need to make to our proof is to the inequalities in Lemma \ref{l:holder-ineq} which we do not carry out here.

Since our estimates, and local existence time are independent of viscosity, we can address the question of convergence in the limit $\nu \to 0$.
\begin{prop}\label{p:u-unu}
Let $u_0 \in \holderspace{k+1}{\alpha}$ be divergence free, and $\ubound$, $T$ be as in Theorem \ref{t:snsexist}. For each $\nu > 0$ we let $u_\nu$ be the solution of the system \snseqns on the time interval $[0,T]$. Making $T$ smaller if necessary, let $u$ be the solution to the Euler equations \eqref{e:euler-evol}--\eqref{e:incompressible} with initial data $u_0$ defined on the time interval $[0,T]$. Then there exists a constant $c = c(k, \alpha, \ubound, L)$ such that for all $t \in [0,T]$ we have
\begin{equation*}
%\hnorm{u(t) - u_\nu(t)}{k}{\alpha} \leqs \frac{c \ubound\sqrt{\nu t}}{L} e^{c \ubound^2 t^2 / L^2}.
\hnorm{u(t) - u_\nu(t)}{k}{\alpha} \leqs \tfrac{c \ubound}{L}\sqrt{\nu t}
\end{equation*}
\end{prop}

At present we are unable to extend the above proposition to domains with boundaries. In this case, possible detachment of the boundary layer creates analytical obstructions to understanding the inviscid limit. We present a proof of Proposition~\ref{p:u-unu} at the end of section \ref{s:snsexist}, and are presently working on extending it to work for domains with boundaries.

\iffalse%{{{ O(t^3/2) theorem
The next theorem we prove shows how the velocity field from the system \snseqns relates to the velocity field of a viscous fluid.
\begin{theorem}\label{t:u-baru}%{{{2 main thm
Let $u_0 \in \holderspace{k+1}{\alpha}$ and $\bar{u}$ be a spatially periodic solution of the Navier-Stokes equations \eqref{e:incompressible}--\eqref{e:ns-evol} with initial data $u_0$. Let $u$ be a solution of the system \snseqns. Then there exists a constant $c = c(k, \alpha, \nu, L, \hnorm{u_0}{k+1}{\alpha})$, and a time $T = T(k, \alpha, \nu, L, \hnorm{u_0}{k+1}{\alpha})$ such that
\begin{equation*}
\hnorm{u(t) - \bar{u}(t)}{k}{\alpha} \leqs \frac{c \hnorm{u_0}{k+1}{\alpha}^2}{L^2} \sqrt{\nu t^3}.
\end{equation*}
for all $t \in [0,T]$.
\end{theorem}

We remark that the super-linear $O(t^{3/2})$ approximation result above is better than a $O(\sqrt{t})$ or $O(t)$ result, which can be obtained elementarily using stationary approximations. We mention again, that we are at present unaware if this estimate is optimal. A proof of Theorem \ref{t:u-baru} can be found in Section \ref{s:snsasym}.
The remainder of this paper is devoted to proving the above results.
\fi

\section{The Weber operator and bounds.}\label{s:wbounds}%{{{1
In this section we define and obtain estimates for the Weber operator which will be central to all subsequent results.\medskip

\begin{defn}%{{{2
We define the Weber operator $\W:\holderspace{k}{\alpha} \times \holderspace{k+1}{\alpha} \to \holderspace{k}{\alpha}$ by
$$ \W( v, \ell) = \lhp\left[ \left(\Imatrix + \gradt \ell\right) v\right] $$
where $\lhp$ is the Leray-Hodge projection \cite{cm} onto divergence free vector fields.
\end{defn}
\begin{rem}\label{r:rangew}%{{{2
The range of $\W$ is \holderspace{k}{\alpha} because multiplication by a $\holderspace{k}{\alpha}$ function is bounded on $\holderspace{k}{\alpha}$, and $\lhp$ is a classical Calderon-Zygmund singular integral operator \cite{stein} which is bounded on H\"older spaces.
\end{rem}
\begin{rem*}%{{{2
In the whole space, or with periodic boundary conditions, the Leray-Hodge projection commutes with derivatives. This is not true for arbitrary domains \cite{constbook}.
\end{rem*}

Formally it seems that $\W(v, \ell)$ should have one less derivative than $\ell$. However we prove below that $\W(v, \ell)$ has as many derivatives as $\ell$. The reason being, when we differentiate $\W(v, \ell)$, we can use `integration by parts' to express the right hand side only in terms of first order derivatives.
\begin{lemma}[Integration by parts]\label{l:ibp}%{{{2
If $u,v \in \holderspace{1}{\alpha}$ then
$$\lhp \left[ \left(\gradt u\right) v\right] = -\lhp \left[ \left(\gradt v\right) u\right]$$
\end{lemma}
\begin{proof}%{{{2
This follows immediately from the identity
\begin{align*}
\left(\gradt u\right) v + \left(\gradt v\right) u = \grad (u \cdot v)
\end{align*}
and the fact that $\lhp$ vanishes on gradients.
\end{proof}

\begin{cor}\label{c:wreg}%{{{2 Range of W
If $k \geqs 1$ and $v, \ell \in \holderspace{k}{\alpha}$ then $\W( v, \ell) \in \holderspace{k}{\alpha}$ and
\begin{equation*}
\hnorm{ \W( v, \ell)}{k}{\alpha} \leqs c \left( 1 + \hnorm{\grad \ell}{k-1}{\alpha} \right) \hnorm{v}{k}{\alpha}.
\end{equation*}
\end{cor}
\begin{proof}%{{{2
Notice first that $\W( v, \ell) \in \holderspace{k-1}{\alpha}$ by Remark \ref{r:rangew}. Now
\begin{align*}
\del_i \W( v, \ell) &= \lhp\left[ (\gradt \del_i \ell) v + \gradt \ell \, \del_i v \right]\\
    &= \lhp\left[ -\gradt v \, \del_i \ell + \gradt \ell \, \del_i v \right].
\end{align*}
Now the right hand side has only first order derivatives of $\ell$ and $v$, hence $\grad \W( v, \ell) \in \holderspace{k-1}{\alpha}$ and the proposition follows.
\end{proof}

\begin{prop}\label{p:wlip}%{{{2
If $k \geqs 1$ and $\ell_1, \ell_2 \in \holderspace{k}{\alpha}$ and $v_1, v_2 \in \holderspace{k}{\alpha}$, are such that
\begin{align*}
\hnorm{\grad\ell_i}{k-1}{\alpha} &\leqs d\\
\text{and}\qquad	\hnorm{v_i}{k}{\alpha} &\leqs \ubound
\end{align*}
for $i = 1,2$, then there exists $c = c( k, d, \alpha)$ such that
\begin{equation}
\label{e:W-lip}	\hnorm{ \W( v_1, \ell_1) - \W(v_2, \ell_2) }{k}{\alpha} \leqs c \left( \tfrac{\ubound}{L} \hnorm{\ell_1 - \ell_2}{k}{\alpha} + \hnorm{v_1 - v_2}{k}{\alpha} \right).
\end{equation}
If $k = 0$, the inequality \eqref{e:W-lip} still holds provided we assume
\begin{align*}
\hnorm{\grad\ell_i}{0}{\alpha} &\leqs d\\
\text{and}\qquad	\hnorm{v_i}{1}{\alpha} &\leqs \ubound
\end{align*}
for $i=1,2$.
\end{prop}%}}}2
\begin{proof}[Proof of Proposition \ref{p:wlip}]%{{{2
The main idea in the proof is to use `integration by parts' to avoid the loss of derivative. By definition of $\W$ we have
\begin{align}
\nonumber \W(v_1, \ell_1) - \W( v_2, \ell_2) &= \lhp\left[ (\Imatrix + \gradt \ell_1) v_1 - (\Imatrix + \gradt \ell_2) v_2 \right]\\
\nonumber    &= \lhp \left[ (\Imatrix + \gradt \ell_1) (v_1 - v_2) + \gradt (\ell_1 - \ell_2) v_2 \right]\\
\label{e:wlip1}    &= \lhp \left[ (\Imatrix + \gradt \ell_1) (v_1 - v_2) - \gradt v_2 (\ell_1 - \ell_2) \right].
\end{align}
Further, differentiating we have
\begin{align}
\nonumber \del_i \left[\W(v_1, \ell_1) - \W( v_2, \ell_2)\right] &= \del_i \lhp \left[ (\Imatrix + \gradt \ell_1) (v_1 - v_2) - \gradt v_2 (\ell_1 - \ell_2) \right]\\
\nonumber &= \lhp \big[ \gradt \del_i\ell_1 (v_1 - v_2) + (\Imatrix + \gradt \ell_1)\del_i(v_2 - v_1) -\\
\nonumber &\qquad - \gradt \del_i v_2 (\ell_1 - \ell_2) - \gradt v_2 \del_i(\ell_2 - \ell_1) \big]\\
\label{e:wlip2} &= \lhp \big[ -\gradt(v_1 - v_2) \del_i\ell_1  + (\Imatrix + \gradt \ell_1)\del_i(v_2 - v_1) +\\
\nonumber &\qquad + \gradt(\ell_1 - \ell_2) \del_i v_2 - \gradt v_2 \del_i(\ell_2 - \ell_1) \big].
\end{align}

Note that we used Lemma \ref{l:ibp} to ensure that the right hand sides of \eqref{e:wlip1} and \eqref{e:wlip2} have only first order derivatives of $\ell$ and $v$. Thus taking the \holderspace{k-1}{\alpha} norms of equations \eqref{e:wlip1} and \eqref{e:wlip2}, and using the fact that multiplication by a \holderspace{k}{\alpha} function and $\lhp$ are bounded on \holderspace{k}{\alpha}, the proposition follows.
\end{proof}%}}}2

\section{Local existence for the stochastic formulation.}\label{s:snsexist}%{{{1
In this section we prove local in time $\holderspace{k}{\alpha}$ existence for the stochastic system \snseqns as stated in Theorem \ref{t:snsexist}. We conclude by proving Proposition \ref{p:u-unu}, showing how the stochastic system \snseqns behaves as $\nu \to 0$. We begin with a few preliminary results.

\begin{lemma}\label{l:holder-ineq}%{{{2
If $k \geqs 1$, then there exists a constant $c = c(k, \alpha)$ such that
\begin{gather*}
\hnorm{f \circ g}{k}{\alpha} \leqs c \hnorm{f}{k}{\alpha}\left(1 + \hnorm{\grad{g}}{k-1}{\alpha}\right)^{k+\alpha}\\
\hnorm{ f\circ g_1 - f \circ g_2}{k}{\alpha} \leqs c \hnorm{\grad f}{k}{\alpha} \left( 1 + \hnorm{\grad g_1}{k-1}{\alpha} + \hnorm{\grad g_2}{k-1}{\alpha} \right)^{k+1} \hnorm{g_1 - g_2}{k}{\alpha}
\end{gather*}
and
\begin{multline*}
\hnorm{f_1 \circ g_1 - f_2 \circ g_2}{k}{\alpha} \leqs c \left( 1 + \hnorm{\grad g_1}{k-1}{\alpha} + \hnorm{\grad g_2}{k-1}{\alpha}\right)^{k+1} \cdot\\
\cdot \left[ \hnorm{f_1 - f_2}{k}{\alpha} + \min\left\{\hnorm{\grad f_1}{k}{\alpha}, \hnorm{\grad f_2}{k}{\alpha} \right\} \hnorm{g_1 - g_2}{k}{\alpha} \right].
\end{multline*}
\end{lemma}%}}}

The proof of Lemma \ref{l:holder-ineq} is elementary and not presented here. We subsequently use the above lemma repeatedly without reference or proof.

\begin{lemma}\label{l:inverse-close}%{{{2
Let  $X_1, X_2 \in \holderspace{k+1}{\alpha}$ be such that
$$\hnorm{\grad X_1 - \Imatrix}{k}{\alpha} \leqs d < 1 \quad\text{and}\quad \hnorm{\grad X_2 - \Imatrix}{k}{\alpha} \leqs d < 1.$$
Let $A_1$ and $A_2$ be the inverse of $X_1$ and $X_2$ respectively. Then there exists a constant $c = c(k, \alpha, d)$ such that
$$\hnorm{A_1 - A_2}{k}{\alpha} \leqs c \hnorm{X_1 - X_2}{k}{\alpha}$$
\end{lemma}

\begin{proof}%{{{2
Let $c = c(k, \alpha, d)$ be a constant that changes from line to line (we use this convention implicitly throughout this paper). Note first $\grad A = (\grad X)\inv \circ A$, and hence by Lemma \ref{l:inverse-bound}
$$\cnorm{\grad A}{0} \leqs \cnorm{(\grad X)\inv}{0} \leqs c.$$
Now using Lemma \ref{l:inverse-bound} to bound $\hnorm{ (\grad X)\inv}{0}{\alpha}$ we have
$$\hnorm{\grad A}{0}{\alpha} = \hnorm{ (\grad X)\inv \circ A}{0}{\alpha} \leqs \hnorm{ (\grad X)\inv}{0}{\alpha} \left( 1 + \cnorm{\grad A}{0} \right) \leqs c$$
When $k \geqs 1$, we again bound $\hnorm{(\grad X)\inv}{k}{\alpha}$ by Lemma \ref{l:inverse-bound}. Taking the \holderspace{k}{\alpha} norm of $(\grad X)\inv \circ A$ we have
$$\hnorm{\grad A}{k}{\alpha} \leqs \hnorm{(\grad X)\inv}{k}{\alpha} \left( 1 + \hnorm{\grad{A}}{k-1}{\alpha} \right)^k.$$
So by induction we can bound $\hnorm{\grad A}{k}{\alpha}$ by a constant $c = c(k, \alpha, d)$. The Lemma now follows immediately from the identity
\begin{align*}
A_1 - A_2 &= \left( A_1 \circ X_2 - \Ifn \right) \circ A_2\\
    &= \left(A_1 \circ X_2 - A_1 \circ X_1 \right) \circ A_2
\end{align*}
and Lemma \ref{l:holder-ineq}.
\end{proof}%}}}
\begin{lemma}\label{l:grad-lambda}%{{{2
Let $u \in C([0,T], \holderspace{k+1}{\alpha})$ and $X$ satisfy the SDE \eqref{e:X-def} with initial data \eqref{e:X-idata}. Let $\lambda = X - \Ifn$ and $\ubound = \sup_t \hnorm{u(t)}{k+1}{\alpha}$. Then there exists $c = c(k, \alpha, \hnorm{u}{k+1}{\alpha})$ such that for short time
\begin{equation*}
\hnorm{\grad \lambda(t)}{k}{\alpha} \leqs \frac{c \ubound t}{L} e^{c\ubound t / L} \qquad\text{and}\qquad\hnorm{\grad \ell(t)}{k}{\alpha} \leqs \frac{c \ubound t}{L} e^{c\ubound t / L}.
\end{equation*}
\end{lemma}

\begin{proof}%{{{2
From equation \eqref{e:X-def} we have
\begin{alignat}{2}
\nonumber	&& X(x,t) &= x + \int_0^t u( X(x,s), s) \,ds + \sqrt{2\nu} B_t\\
\label{e:gradx-def}	&\implies\quad&	\grad X(t) &= I + \int_0^t (\grad u)\circ X \cdot \grad X.
\end{alignat}
Taking the \cspace{0} norm of equation \eqref{e:gradx-def} and using Gronwall's Lemma we have
\begin{equation*}
\cnorm{\grad \lambda(t)}{0} = \cnorm{ \grad X(t) - I}{0} \leqs e^{\ubound t/L} - 1.
\end{equation*}
Now taking the \holderspace{k}{\alpha} norm in equation \eqref{e:gradx-def} we have
$$ \hnorm{ \grad \lambda(t) }{k}{\alpha} \leqs  c \int_0^t \hnorm{\grad u}{k}{\alpha} \left(1 + \hnorm{\grad \lambda}{k-1}{\alpha} \right)^k \left( 1 + \hnorm{ \grad \lambda}{k}{\alpha}\right).$$
The bound for $\hnorm{\grad \lambda}{k}{\alpha}$ now follows from the previous two inequalities, induction and Gronwall's Lemma. The bound for $\hnorm{\grad \ell}{k}{\alpha}$ then follows from Lemma \ref{l:inverse-close}.

We draw attention to the fact that the above argument can only bound $\grad \lambda$, and not $\lambda$. Fortunately, our results only rely on a bound of $\grad \lambda$.
\end{proof}

\begin{lemma}\label{l:X-gronwall}%{{{2
Let $u, \bar{u} \in C( [0,T], \holderspace{k+1}{\alpha})$ be such that
$$ \sup_{0 \leqs t \leqs T} \hnorm{u(t)}{k+1}{\alpha} \leqs \ubound \quad\text{and}\quad \sup_{0 \leqs t \leqs T} \hnorm{\bar{u}(t)}{k+1}{\alpha} \leqs \ubound.$$
Let $X, \bar{X}$ be solutions of the SDE \eqref{e:X-def}--\eqref{e:X-idata} with drift $u$ and $\bar{u}$ respectively, and let $A$ and $\bar{A}$ be the spatial inverse of $X$ and $\bar{X}$ respectively. Then there exists $c = c(k, \alpha, \ubound)$ and a time $T=T(k, \alpha, \ubound)$ such that
\begin{gather}
\label{e:X-Xbar} \hnorm{X(t) - \bar{X}(t)}{k}{\alpha} \leqs c e^{c \ubound t/L} \int_0^t \hnorm{u - \bar{u}}{k}{\alpha}\\
\label{e:A-Abar} \hnorm{A(t) - \bar{A}(t)}{k}{\alpha} \leqs c e^{c \ubound t/L} \int_0^t \hnorm{u - \bar{u}}{k}{\alpha}
\end{gather}
for all $0 \leqs t \leqs T'$.
\end{lemma}
\begin{proof}%{{{2
We first use Lemma \ref{l:grad-lambda} to bound $\hnorm{\grad X - \Imatrix}{k}{\alpha}$ and $\hnorm{\grad \bar{X} - \Imatrix}{k}{\alpha}$ for short time $T'$. Now
\begin{align*}
&&	X(t) - \bar{X}(t) &= \int_0^t u\circ X - \bar{u}\circ \bar{X}\\
&\implies&	\hnorm{X(t) - \bar{X}(t)}{k}{\alpha} &\leqs \int_0^t \hnorm{u\circ X - \bar{u}\circ \bar{X}}{k}{\alpha}\\
&&	&\leqs c\int_0^t \left( \hnorm{u - \bar{u}}{k}{\alpha} + \frac{\ubound}{L}\hnorm{X - \bar{X}}{k}{\alpha}\right)
\end{align*}
and inequality \eqref{e:X-Xbar} follows by applying Gronwall's Lemma. Inequality \eqref{e:A-Abar} follows immediately from \eqref{e:X-Xbar} and Lemma \ref{l:inverse-close}.
\end{proof}%}}}

We now provide the proof of Theorem \ref{t:snsexist}. We reproduce the statement here for convenience.
\begin{theorem*}[\ref{t:snsexist}]%{{{2 Restate t:snsexist
Let $k \geqs 1$ and $u_0 \in \holderspace{k+1}{\alpha}$ be divergence free. There exists a time $T = T(k, \alpha, L, \hnorm{u_0}{k+1}{\alpha})$, but independent of viscosity, and a pair of functions $\lambda, u \in C([0,T], \holderspace{k+1}{\alpha})$ such that $u$ and $X = \Ifn + \lambda$ satisfy the system \snseqns. Further $\exists \ubound = \ubound(k, \alpha, L, \hnorm{u_0}{k+1}{\alpha})$ such that $t \in [0,T] \implies \hnorm{u(t)}{k+1}{\alpha} \leqs \ubound$.
\end{theorem*}
%\begin{proof}[Proof of Theorem \ref{t:snsexist}]%{{{2
\begin{proof}
Let $\ubound$ be a large constant, and $T$ a small time, both of which will be specified later. Define as before $\uball$ and $\lball$ by
\begin{align*}
\uball &= \left\{ u \in C([0,T], \holderspace{k+1}{\alpha}) \;\big|\; \hnorm{u(t)}{k+1}{\alpha} \leqs \ubound, \; \divergence u = 0 \text{ and } u(0) = u_0 \right\}\\
\text{and}\quad \lball &= \left\{ \ell \in C([0,T], \holderspace{k+1}{\alpha}) \;\big|\; \hnorm{\grad \ell(t)}{k}{\alpha} \leqs \tfrac{1}{2} \; \forall t \in [0,T] \text{ and } \ell(\cdot, 0) = 0 \right\}.
%\lball = \left\{ \ell \in \holderspace{k+1}{\alpha}(\I \times [0,T], \I) \;\big|\; \hnorm{\grad \ell(t)}{k}{\alpha} \leqs \tfrac{1}{2} \; \forall t \in [0,T] \text{ and } \ell(\cdot, 0) = 0 \right\}.
\end{align*}
We clarify that the functions $u$ and $\ell$ are required to be spatially $\holderspace{k+1}{\alpha}$, and need only be continuous in time.

Now given $u \in \uball$ we define $X_u$ to be the solution of equation \eqref{e:X-def} with initial data \eqref{e:X-idata} and $\lambda_u = X_u - I$ be the Eulerian displacement. We define $A_u$ by equation \eqref{e:A-def} and let $\ell_u = A_u - I$ be the Lagrangian displacement. Finally we define $W: \uball \to \uball$ by
\begin{equation*}
W(u) = \E\W(u_0 \circ A_u, \ell_u).
\end{equation*}
We aim to show that $W: \uball \to \uball$ is Lipschitz in the weaker norm
\begin{equation*}
\norm{u}_\uball = \sup_{0 \leqs t \leqs T} \hnorm{u(t)}{k}{\alpha}
\end{equation*}
and when $T$ is small enough, we will show that $W$ is a contraction mapping.\smallskip

Let $c$ be a constant that changes from line to line. By Corollary \ref{c:wreg} we have
\begin{align}
\nonumber \hnorm{W(u)}{k+1}{\alpha} &\leqs c \E\left[ \left( 1 + \hnorm{\grad \ell_u}{k}{\alpha} \right) \hnorm{u_0 \circ A_u}{k+1}{\alpha}\right]\\
\label{e:wubd}    & \leqs c \hnorm{u_0}{k+1}{\alpha} \sup_\Omega \left( 1 + \hnorm{\grad\ell_u}{k}{\alpha} \right)^{k+2}.
\end{align}
Here $\Omega$ is the probability space on which our processes are defined. We remark that Lemma \ref{l:grad-lambda} gives us a bound on $\hnorm{\grad \ell_u}{k}{\alpha}$. A bound on $\E \hnorm{\grad \ell_u}{k}{\alpha}$ instead would not have been enough.

Now we choose $\ubound = c(\frac{3}{2})^{k+2}\hnorm{u_0}{k+1}{\alpha}$, and then apply Lemma \ref{l:grad-lambda} to choose $T$ small enough to ensure $\ell_u,\lambda_u \in \lball$. Now inequality \ref{e:wubd} ensures that $W(u) \in \uball$. Now if $u,\bar{u} \in \uball$, Lemma \ref{l:X-gronwall} guarantees
\begin{equation*}
\hnorm{\ell_u(t) - \ell_{\bar{u}}(t)}{k}{\alpha} \leqs c e^{c\ubound t/L} \int_0^t \hnorm{u - \bar{u}}{k}{\alpha}.
\end{equation*}
Thus applying Proposition \ref{p:wlip} we have
\begin{align*}
\hnorm{W(u)(t) - W(\bar{u})(t)}{k}{\alpha} &\leqs c \left( \tfrac{\ubound}{L} \hnorm{ \ell_u(t) - \ell_{\bar{u}}(t)}{k}{\alpha} + \hnorm{u_0 \circ A_u(t) - u_0 \circ A_{\bar{u}}(t)}{k}{\alpha} \right)\\
    &\leqs \frac{c \ubound}{L} \hnorm{ \ell_u(t) - \ell_{\bar{u}}(t)}{k}{\alpha}\\
    &\leqs \frac{c \ubound}{L} e^{c \ubound t/L} \int_0^t \hnorm{u - \bar{u}}{k}{\alpha}.
\end{align*}
So choosing $T=T(k, \alpha, L, \ubound)$ small enough we can ensure $W$ is a contraction.\medskip

The existence of a fixed point of $W$ now follows by successive iteration. We define $u_{n+1} = W(u_n)$. The sequence $(u_n)$ converges strongly with respect to the $\holderspace{k}{\alpha}$ norm. Since $\uball$ is closed and convex, and the sequence $(u_n)$ is uniformly bounded in the $\holderspace{k+1}{\alpha}$ norm, it must have a weak limit $u \in \uball$. Finally since $W$ is continuous with respect to the weaker $\holderspace{k}{\alpha}$ norm, the limit must be a fixed point of $W$, and hence a solution to the system \snseqns.
\end{proof}
%}}}2
% TODO: State and motivate the other model.
% \begin{rem}%{{{2 Coefficients in the evolution equation of v
% We note that the commutator coefficients (equation \eqref{e:cdef}) only depend on spatial gradients of $A$, and not on $A$ itself, or on time derivatives of $A$. The process $A$ can be arbitarily large with non-zero probability and is almost surely not differentiable in time, however the spatial gradient is well behaved (as shown below). Thus local existence, and the parabolic regularity estimates for equation \eqref{e:v-def} hold uniformly in $\Omega$, and the webber operator is well defined (and proposition \ref{p:wns-lip} applies).
% \end{rem}

We conclude by proving the vanishing viscosity behavior stated in Proposition~\ref{p:u-unu}. We reproduce the statement here for convenience.
\begin{prop*}[\ref{p:u-unu}]%{{{ Restate u - u_nu prop
Let $u_0 \in \holderspace{k+1}{\alpha}$ be divergence free, and $\ubound$, $T$ be as in Theorem \ref{t:snsexist}. For each $\nu > 0$ we let $u_\nu$ be the solution of the system \snseqns on the time interval $[0,T]$. Making $T$ smaller if necessary, let $u$ be the solution to the Euler equations \eqref{e:euler-evol}--\eqref{e:incompressible} with initial data $u_0$ defined on the time interval $[0,T]$. Then there exists a constant $c = c(k, \alpha, \ubound, L)$ such that for all $t \in [0,T]$ we have
\begin{equation*}
%\hnorm{u(t) - u_\nu(t)}{k}{\alpha} \leqs \frac{c \ubound\sqrt{\nu t}}{L} e^{c \ubound^2 t^2 / L^2}.
\hnorm{u(t) - u_\nu(t)}{k}{\alpha} \leqs \tfrac{c \ubound}{L}\sqrt{\nu t}
\end{equation*}
\end{prop*}
\begin{proof}%{{{2
We use a subscript of $\nu$ to denote quantities associated to the solution of viscous problem \snseqns, and unsubscripted letters to denote the corresponding quantities associated to the solution of the Eulerian-Lagrangian formulation of the Euler equations \eqref{e:ele-adef}--\eqref{e:ele-webber}. We use the same notation as in the proof of Theorem~\ref{t:snsexist}.

Now from the proof of Theorem~\ref{t:snsexist} we know that for short time $\ell_\nu, \ell \in \lball$. Using Lemma \ref{l:inverse-close} and making $T$ smaller if necessary, we can ensure $\lambda_\nu, \lambda \in \lball$. We begin by estimating $\E\hnorm{\lambda_\nu - \lambda}{k}{\alpha}$:
\begin{alignat*}{2}
&&	\lambda_\nu(t) - \lambda(t) &= \int_0^t \left[ u_\nu \circ X_\nu - u \circ X \right] + \sqrt{2\nu}B_t\\
&\implies\quad&	\hnorm{\lambda_\nu(t) - \lambda(t)}{k}{\alpha} &\leqs c \left( \int_0^t \left[ \hnorm{u_\nu - u}{k}{\alpha} + \tfrac{\ubound}{L} \hnorm{ \lambda_\nu - \lambda}{k}{\alpha}  \right] + \sqrt{\nu} |B_t| \right)
\end{alignat*}
and so by Gronwall's lemma
\begin{equation*}
\hnorm{\lambda_\nu(t) - \lambda(t)}{k}{\alpha} \leqs c \left( \sqrt{\nu} |B_t| + \int_0^t \hnorm{u_\nu - u}{k}{\alpha} \right) e^{c \ubound t / L}.
\end{equation*}
Using Lemma \ref{l:inverse-close} and taking expected values gives
\begin{equation}
\label{e:lnu-l} \E \hnorm{\ell_\nu(t) - \ell(t)}{k}{\alpha} \leqs c \left( \sqrt{\nu t} + \int_0^t \hnorm{u_\nu - u}{k}{\alpha} \right) e^{c \ubound t / L}.
\end{equation}

To estimate the difference $u_\nu - u$, we use \eqref{e:u-def}, and \eqref{e:ele-webber} to obtain
\begin{alignat*}{2}
&&	u_\nu - u &= \E\W( u_0 \circ A_\nu, \ell_\nu) - \W(u_0 \circ A, \ell)\\
&\implies\quad& \hnorm{ u_\nu - u }{k}{\alpha} &\leqs c \E\left( \tfrac{\ubound}{L} \hnorm{\ell_\nu - \ell}{k}{\alpha} + \hnorm{u_0\circ A_\nu - u_0 \circ A}{k}{\alpha} \right)\\
&&	&\leqs \tfrac{c \ubound}{L} \E\hnorm{\ell_\nu - \ell}{k}{\alpha}\\
&\implies&\hnorm{ u_\nu(t) - u(t) }{k}{\alpha}	&\leqs \tfrac{c \ubound}{L} e^{c \ubound t/L} \left( \sqrt{\nu t} + \int_0^t \hnorm{u_\nu - u}{k}{\alpha} \right)
\end{alignat*}
and the theorem follows from Gronwall's lemma.
\end{proof}
\section{Local existence for the Navier-Stokes equations.}\label{s:nsexist}%{{{1
Proposition \ref{p:wlip}, along with Peter Constantin's diffusive Lagrangian formulation \cite{elns} immediately gives us a local in time \holderspace{k}{\alpha} existence and uniqueness result for the Navier-Stokes equations using classical PDE methods. We conclude this paper by presenting the proof in this section.

\begin{defn}\label{d:wns}%{{{2
Let $k \geqs 2$ and $T>0$. We define $\lballT{k}{\alpha}{T}$ by
$$\lballT{k}{\alpha}{T} = \left\{ \ell \in \holderspace{k}{\alpha}(\I \times [0,T], \I) \;\big|\; \hnorm{\grad \ell(t)}{k-1}{\alpha} \leqs \tfrac{1}{2} \; \forall t \in [0,T] \text{ and } \ell(\cdot, 0) = 0 \right\}.$$
Given $\ell \in \lballT{k}{\alpha}{T}$, and $u \in \holderspace{k}{\alpha}( \I\times[0,T], \I)$ divergence free we define the \textit{virtual velocity} $v = v_{u,\ell}$ to be the unique solution of the linear parabolic equation
\begin{equation}
\label{e:v-def}	\D{v_\beta} = 2 \nu C^i_{j,\beta} \del_j v_i
\end{equation}
with initial data
\begin{equation}
\label{e:vidata}	v(x,0) = u(x,0)
\end{equation}
where
\begin{equation}
\label{e:cdef}	C^p_{j,i} = (\Imatrix + \grad \ell)\inv_{ki} \del_k \del_j \ell_p
\end{equation}
are the commutator coefficients.

Finally we define the operator $\Wns:\holderspace{k}{\alpha}(\I \times[0,T]) \times \lballT{k}{\alpha}{T} \to \holderspace{k}{\alpha}(\I \times [0,T])$ by
\begin{equation}
\label{e:wns-def} \Wns( u, \ell) = \W ( v_{u, \ell}, \ell).
\end{equation}
\end{defn}
\begin{rem*}%{{{2
We clarify that by $\ell \in \holderspace{k}{\alpha}(\I \times [0,T], \I)$ we only impose a $\holderspace{k}{\alpha}$ \textit{spatial} regularity restriction. We do not assume anything about time regularity. This will be the case for the remainder of this section.
\end{rem*}
\begin{rem*}%{{{2
Observe that $\hnorm{\grad\ell}{k-1}{\alpha} \leqs \frac{1}{2}$ guarantees that the matrix $\Imatrix + \gradt\ell$ in equation \eqref{e:cdef} is invertible. Further note that all coefficients in equation \eqref{e:v-def} are of class \holderspace{k}{\alpha} and hence by parabolic regularity~\cite{krylov}, $v \in \holderspace{k}{\alpha}$.
\end{rem*}

\begin{lemma}\label{l:inverse-bound}%{{{2
Let $X$ be a Banach algebra. If $x \in X$ is such that $\norm{x} \leqs \rho < 1$ then $1 + x$ is invertible and $\norm{(1 + x)\inv} \leqs \frac{1}{1 - \rho}$. Further if in addition $\norm{y} \leqs \rho$ then
$$\norm{(1+x)\inv - (1+y)\inv} \leqs \frac{1}{(1 - \rho)^2} \norm{x-y}$$
\end{lemma}
\begin{proof}%{{{2
The first part of the Lemma follows immediately from the identity $(1+x)\inv = \sum (-x)^n$. The second part follows from the first part and the identity
\begin{equation*}
(1+x)\inv - (1+y)\inv = (1+x)\inv (y-x) (1+y)\inv.\qedhere
\end{equation*}
\end{proof}%}}}

We generally use Lemma \ref{l:inverse-bound} when $X$ is the space of \holderspace{k}{\alpha} periodic matrices. We finally prove that the Weber operator $\Wns$ is Lipschitz, which will quickly give us the existence theorem.

\begin{prop}\label{p:wns-lip}%{{{2 W bounds
If $\ell,\bar{\ell} \in \lballT{k}{\alpha}{T}$ and $u, \bar{u} \in \holderspace{k}{\alpha}$, are such that
$$\sup_{0 \leqs t \leqs T} \hnorm{u(t)}{k}{\alpha} \leqs \ubound \quad \text{and} \quad \sup_{0 \leqs t \leqs T} \hnorm{\bar{u}(t)}{k}{\alpha} \leqs \ubound$$
then there exists $c = c( k, \alpha, L, \nu, \ubound)$ and $T' = T'(k, \alpha, L, \nu, \ubound)$ such that
\begin{multline*}
\hnorm{\Wns(u, \ell)(t) - \Wns(\bar{u}, \bar{\ell})(t)}{k}{\alpha} \leqs c \hnorm{u(0) - \bar{u}(0)}{k}{\alpha} + \\
+ \frac{c \ubound}{L} \left[ \left( 1 + \frac{\nu t}{L^2} \right) \hnorm{\ell(t) - \bar{\ell}(t)}{k}{\alpha} + t \hnorm{u(t) - \bar{u}(t)}{k-2}{\alpha} \right]
\end{multline*}
for all $0 \leqs t \leqs T'$.
\end{prop}

\begin{proof}%{{{2
Let $v$ and $\bar{v}$ be the virtual velocities associated to $u, \ell$ and $\bar{u}, \bar{\ell}$ respectively. Let $C$ and $\bar{C}$ be the commutator coefficients associated to $\ell$ and $\bar{\ell}$ respectively. Since equation \eqref{e:v-def} is a linear parabolic equation with \holderspace{k}{\alpha} coefficients, standard regularity theory \cite{krylov} ensures that there exists $T' = T'( \nu, \ubound)$ such that
\begin{equation*}
\sup_{0 \leqs t \leqs T'} \hnorm{v(t)}{k}{\alpha} \leqs 2\ubound \quad\text{and}\quad	\sup_{0 \leqs t \leqs T'} \hnorm{\bar{v}(t)}{k}{\alpha} \leqs 2\ubound.
\end{equation*}
Hence by proposition \ref{p:wlip} we have
\begin{align}
\nonumber \hnorm{\Wns(u, \ell) - \Wns( \bar{u}, \bar{\ell})}{k}{\alpha} &= \hnorm{\W( v, \ell) - \W( \bar{v}, \bar{\ell})}{k}{\alpha}\\
\label{e:wtmp}    &\leqs c \left( \tfrac{\ubound}{L} \hnorm{ \ell - \bar{\ell}}{k}{\alpha} + \hnorm{v - \bar{v}}{k}{\alpha} \right)
\end{align}

Now let $\tilde{v} = v - \bar{v}$. The evolution equation of $\tilde{v}$ is given by
$$\D[u]{\tilde{v}_\beta} - 2 \nu C^i_{j,\beta} \del_j \tilde{v}_i = 2 \nu (\bar{C}^i_{j,\beta} - C^i_{j,\beta} ) \del_j \bar{v}_i + ((\bar{u} - u) \cdot \grad ) \bar{v}_\beta$$
with initial data
$$\tilde{v}(x,0) = u(x,0) - \bar{u}(x,0).$$
We estimate the $\holderspace{k-2}{\alpha}$ norm of the right hand side. Let $c$ be some constant which changes from line to line. By definition,
\begin{align*}
\bar{C}^k_j - C^k_j &=  (\Imatrix + \gradt \bar{\ell})^{-1} \grad \del_j \bar{\ell}_k - (\Imatrix + \gradt \ell)^{-1} \grad \del_j \ell_k \\
    &= \left[ (\Imatrix + \gradt \bar{\ell})^{-1} - (\Imatrix + \gradt \ell)\inv \right] \grad \del_j \bar{\ell}_k + (\Imatrix + \gradt \ell)\inv \left[ \grad \del_j \bar{\ell}_k - \grad \del_j \ell_k \right].
\end{align*}
Note that by Lemma \ref{l:inverse-bound} we can bound $\hnorm{ (\Imatrix + \gradt \ell)\inv}{k-1}{\alpha}$ and $\hnorm{ (\Imatrix + \gradt \bar{\ell})\inv}{k-1}{\alpha}$. Further, by Lemma \ref{l:inverse-bound} again we have
$$
\hnorm{(\Imatrix + \gradt \bar{\ell})\inv - (\Imatrix + \gradt \ell)\inv}{k-1}{\alpha} \leqs c \hnorm{\grad\ell - \grad\bar{\ell}}{k-1}{\alpha}.
$$
Combining these estimates we have
\begin{align*}
\hnorm{C - \bar{C}}{k-2}{\alpha} \leqs \frac{c}{L} \hnorm{\grad\ell - \grad\bar{\ell}}{k-1}{\alpha}.
\end{align*}
Finally note
\begin{equation*}
\hnorm{( ( \bar{u} - u) \cdot \grad) \bar{v}}{k-2}{\alpha} \leqs \frac{c \ubound}{L} \hnorm{u - \bar{u}}{k-2}{\alpha}
\end{equation*}
Thus by parabolic regularity \cite{krylov},
\begin{multline}
\label{e:tildevbd} \hnorm{\tilde{v}(t)}{k}{\alpha} \leqs \frac{ c \ubound t}{L} \left( \frac{\nu}{L} \hnorm{\grad\ell(t) - \grad\bar{\ell}(t)}{k-1}{\alpha} + \hnorm{u(t) - \bar{u}(t)}{k-2}{\alpha} \right) +\\
+ \hnorm{u(0) - \bar{u}(0)}{k}{\alpha}
\end{multline}
and substituting equation \eqref{e:tildevbd} in \eqref{e:wtmp}, the proposition follows.
\end{proof}

\begin{theorem}\label{t:lens}%{{{2 Local existence for Navier Stokes.
Let $k\geqs 2$ and $u_0 \in \holderspace{k}{\alpha}(\I,\I)$ be divergence free. Then there exists $T = T(k, \alpha, L, \nu, \hnorm{u_0}{k}{\alpha})$ and $u \in \holderspace{k}{\alpha}(\I \times [0,T], \I)$ which is a solution of the Navier-Stokes equations with initial data $u_0$.
\end{theorem}

\begin{proof}%{{{2
Let $\ubound > \hnorm{u_0}{k}{\alpha}$. We define the set $\uball$ by
$$\uball = \left\{ u \in C( [0,T], \holderspace{k}{\alpha}) \;\big|\; \hnorm{u(t)}{k}{\alpha} \leqs \ubound, \; \divergence u = 0, \text{ and } u(0) = u_0 \right\}.$$
Given $u \in \uball$, let $\ell_u$ to be the unique solution of the equation
\begin{equation*}
\D{\ell_u} + u = 0
\end{equation*}
with initial data
\begin{equation*}
\ell_u(x,0) = 0.
\end{equation*}
Our aim is to produce $u \in \uball$ such that $u = \Wns( u, \ell_u)$, which from \cite{elns} we know must be a solution to the Navier-Stokes equations.

We define the map $W$ by
\begin{equation*}
W(u) = \Wns(u, \ell_u).
\end{equation*}
If $\uball$ is endowed with the strong norm
\begin{equation*}
\norm{u}_\uball = \sup_{0 \leqs t \leqs T} \hnorm{u(t)}{k}{\alpha}
\end{equation*}
we will show as before that for sufficiently small $T$, $W$ maps the $\uball$ into itself. Finally we will show that $W$ is a contraction under a weaker norm, producing the desired fixed point.\smallskip

First note that by parabolic regularity \cite{krylov}, we have
\begin{equation*}
\hnorm{ \ell_u(t) }{k}{\alpha} \leqs c \ubound t.
\end{equation*}
The constant $c$ of course depends on $\ubound$, but we retain the $\ubound$ on the right for dimensional correctness. Thus choosing $T$ small will guarantee $\ell_u \in \lball^{k,\alpha}_T$.

Let $v = v_{\ell, u}$ be the virtual velocity defined by equation \eqref{e:v-def}, with initial data $u_0$. Standard parabolic estimates \cite{krylov}, (and the fact that $\ell \in \lball^{k,\alpha}_T$), show
\begin{equation*}
\hnorm{v(t) - u_0}{k}{\alpha} \leqs \frac{c \ubound^2}{L} t
\end{equation*}

Now by definition,
\begin{align*}
W(u) &= \lhp \left[ \left( \Imatrix + \gradt \ell_u \right) v \right]\\
    &= \lhp \left[  \left( \Imatrix + \gradt \ell_u \right) \left( v - u_0 \right) + \left( \Imatrix + \gradt \ell_u \right) u_0 \right]\\
    &= \lhp[ u_0 ] + \lhp\left[ \left(\Imatrix + \gradt \ell_u \right) (v - u_0) \right] + \lhp[ (\gradt \ell_u) u_0 ]
\end{align*}
Since $u_0$ is divergence free, $\lhp(u_0) = u_0$. Using Corollary \ref{c:wreg}, the preceding two estimates for $\ell_u$ and $v - u_0$, we obtain
\begin{align*}
\hnorm{W(u)(t)}{k}{\alpha} \leqs \hnorm{u_0}{k}{\alpha} + c \frac{\ubound^2}{L} t.
\end{align*}
Thus choosing $T < \frac{L}{c \ubound^2} (\ubound - \hnorm{u_0}{k}{\alpha})$, we can ensure that $W$ maps $\uball$ into itself.\smallskip

To see that $W$ has a fixed point, let $u,\bar{u} \in \uball$ and define $\tilde{\ell} = \ell_u - \ell_{\bar{u}}$. The evolution of $\tilde{\ell}$ is governed by
\begin{equation*}
\D{\tilde{\ell}} = \left( (\bar{u} - u) \cdot \grad \right) \ell_{\bar{u}} + \bar{u} - u
\end{equation*}
and parabolic regularity \cite{krylov} immediately gives
\begin{equation*}
\hnorm{\tilde{\ell}(t)}{k}{\alpha} \leqs c t \hnorm{u(t) - \bar{u}(t)}{k-2}{\alpha}
\end{equation*}
Combining this with Proposition \ref{p:wns-lip}, we have
\begin{equation*}
\sup_{0 \leqs t \leqs T'} \hnorm{W(u)(t) - W(\bar{u})(t)}{k}{\alpha} \leqs \frac{c \ubound T'}{L} \sup_{0 \leqs t \leqs T'} \hnorm{u(t) - \bar{u}(t)}{k-2}{\alpha}.
\end{equation*}
Thus if $T'$ is chosen to be smaller than $\frac{L}{c\ubound}$ then $W: \uball \to \uball$ is a contraction mapping and has a unique fixed point concluding the proof.
\end{proof}%}}}2

\begin{rem*}%{{{2
The above estimates along with the active vector formulation of the Euler equations \cite{ele} can be used to prove a \holderspace{k}{\alpha} local existence and uniqueness theorem for the Euler equations. Since a similar proof of this result can be found in the original paper \cite{ele} by P. Constantin, we do not present it here.
\end{rem*}

\section*{Acknowledgment}%{{{1
I would like to thank Peter Constantin for his encouragement, support and many helpful discussions. I would also like to thank Hongjie Dong and Tu Nguyen for carefully reading this paper, and pointing out an error in the original proof of Theorem~\ref{t:snsexist}.
%}}}
% \appendix
% \section{Basic inequalities on H\"older spaces}%{{{1 Stuff for me only
% \begin{lemma}
% There exists a constant $c$ such that
% \begin{gather*}
% \hsnorm{f \circ g}{\alpha} \leqs c \cnorm{\grad f}{0} \hsnorm{g}{\alpha}\\
% \hnorm{f \circ g}{0}{\alpha} \leqs c \left( \cnorm{f}{0} + \cnorm{\grad f}{0} \hsnorm{g}{\alpha} \right) \leqs c \cnorm{f}{1} \left(1 +\frac{\hsnorm{g}{\alpha}}{L} \right)\\
% \hnorm{f \circ g}{0}{\alpha} \leqs c \hnorm{f}{0}{\alpha} \left(1 + \cnorm{ \grad g}{0}^\alpha \right)\\
% \cnorm{f \circ g}{k} \leqs c \cnorm{f}{k}(1 + \cnorm{\grad{g}}{k-1})^k \\
% \hnorm{f \circ g}{k}{\alpha} \leqs c \hnorm{f}{k}{\alpha}\left(1 + \hnorm{\grad{g}}{k-1}{\alpha}\right)^k \\
% \hnorm{ f\circ g_1 - f \circ g_2}{2}{\alpha} \leqs c \cnorm{\grad f}{3} \left( 1 + \cnorm{\grad g_1}{2} + \cnorm{\grad g_2}{2} \right)^3 \hnorm{g_1 - g_2}{2}{\alpha}\\
% \begin{split}
% \hnorm{f_1 \circ g_1 - f_2 \circ g_2}{2}{\alpha} &\leqs c \left( 1 + \cnorm{\grad g_1}{2} + \cnorm{\grad g_2}{2}\right)^3\\
% &\quad\quad\left[ \hnorm{f_1 - f_2}{2}{\alpha} + \cnorm{\grad f}{3} \hnorm{g_1 - g_2}{2}{\alpha} \right]
% \end{split}
% \end{gather*}
% \end{lemma}
% %}}}1
%{{{1 Bibliography
%}}}1
\end{document}